\documentclass[12pt]{article}
\usepackage{fancyhdr,graphicx}
\usepackage{amsfonts}
\usepackage{amssymb}
\usepackage{amsthm}
\usepackage{newlfont}
\usepackage{amsmath}
\usepackage{cases}

\newtheorem{thm}{Theorem}[section]
\newtheorem{Lemma}[thm]{Lemma}
\newtheorem{cor}[thm]{Corollary}

\DeclareMathOperator\dg{deg}
\DeclareMathOperator\dist{dist}

\numberwithin{equation}{section}
\usepackage{latexsym, bm}

\begin{document}
\title{\bf{Connectivity and Minimal Distance Spectral Radius of Graphs }}
\author{Xiaoling Zhang$^{ab}$
\thanks{E-mail addresses: zhangxling04@lzu.cn, cgodsil@uwaterloo.ca.}
~\thanks{The author is currently a visiting Ph.D. student at the Department
 of Combinatorics and Optimisation in University of Waterloo from September 2008 to September 2009.
She is partially supported by NSFC grant no. 10831001.},~~Chris Godsil$^{a}$\footnotemark[1]\\
\small{$^{a}$Combinatorics and Optimisation
    University of Waterloo,}
   \\\small{ Waterloo, Ontario,
    Canada, N2L 3G1}
    \\\small{$^{b}$School of Mathematics and Statistics, Lanzhou University,}
\\\small{ Lanzhou, Gansu 730000, P. R. China}}

\date{}


\maketitle

\begin{abstract}
In this paper, we study how the distance spectral radius behaves
when the graph is perturbed by grafting edges.  As  applications,
we also determine the graph with $k$ cut vertices (respectively,
$k$ cut edges) with the minimal distance spectral radius.

 ~~~~~~~~~~~


~~~~~~~~~~~


\end{abstract}

\section{Introduction}

The distance matrix of a graph, while not as common as the more
familiar adjacency matrix, has nevertheless come up in several
different areas, including communication network design \cite{d},
graph embedding theory \cite{c,e,f}, molecular stability \cite{g,h},
and network flow algorithms \cite{a,b}. So it is interesting to
study the spectra of these matrices. In this paper, we study the
largest eigenvalue of the distance matrix of a graph.

Throughout this paper, we will assume that $G$ is a simple,
connected graph of order $n$, that is, with $n$ vertices.  Let $G$
be a connected graph with vertex set $\{1,\ldots, n\}$. The distance
between vertices $i$ and $j$ of $G$, denoted by $\dist(i, j)$, is
defined to be the length (i.e., the number of edges) of the shortest
path from $i$ to $j$ \cite{am}. The distance matrix of $G$, denoted
by $D(G)$ is the $n \times n$ matrix with its $(i,j)$-entry equal to
$\dist(i, j)$, $i, j = 1, 2, \ldots, n$. Note that $\dist(i, i) =
0$, $i = 1, 2, \ldots, n$. The distance eigenvalue of largest
magnitude is called the distance spectral radius, and is denoted by
$\Lambda_1$. Balaban et al. \cite{al} proposed the use of
$\Lambda_1$ as a structure-descriptor, and it was successfully used
to make inferences about the extent of branching and boiling points
of alkanes\cite{al,an}.

In this paper, we determine the graph with $k$ cut vertices
(respectively, $k$ cut edges) which has the minimal distance
spectral radius.

\section*{Main results}

Let $G$ be a connected graph. Let $\dg(v)$ (or $\dg_G(u)$) denote the degree of
the vertex $v$ in $G$.  We define a pendant path of $G$ to be a
walk  $v_{0}v_{1}\cdots v_{s}$ $(s \geqslant 1)$ such that the
vertices $v_{0}, v_{1}, \ldots, v_{s}$ are distinct, $\dg(v_{0}) > 2$,
$\dg(v_{s}) = 1$, and $\dg(v_{i})= 2$, whenever $0 < i < s$. And
$v_{0}$, $s$  are called the root and  the length of the pendant
path, respectively.

We give a generalization of Theorem 3.5 in \cite{m}.

\begin{thm} \label{1}
Let $u$ and $v$ be two adjacent vertices of a connected graph $G$
and for positive integers $k$ and $l$, let $G_{k, l}$ denote the
graph obtained from $G$ by adding  paths of length $k$ at $u$ and
length $l$ at $v$. If $k > l\geqslant 1$, then $\Lambda_1(G_{k, l})<
\Lambda_1(G_{k+1, l-1})$; if $k=l\geqslant 1$, then $\Lambda_1(G_{k,
l})< \Lambda_1(G_{k+1, l-1})$ or $\Lambda_1(G_{k, l})<
\Lambda_1(G_{k-1, l+1})$.
\end{thm}

\begin{figure}[!htbp]
\begin{center}
    \includegraphics{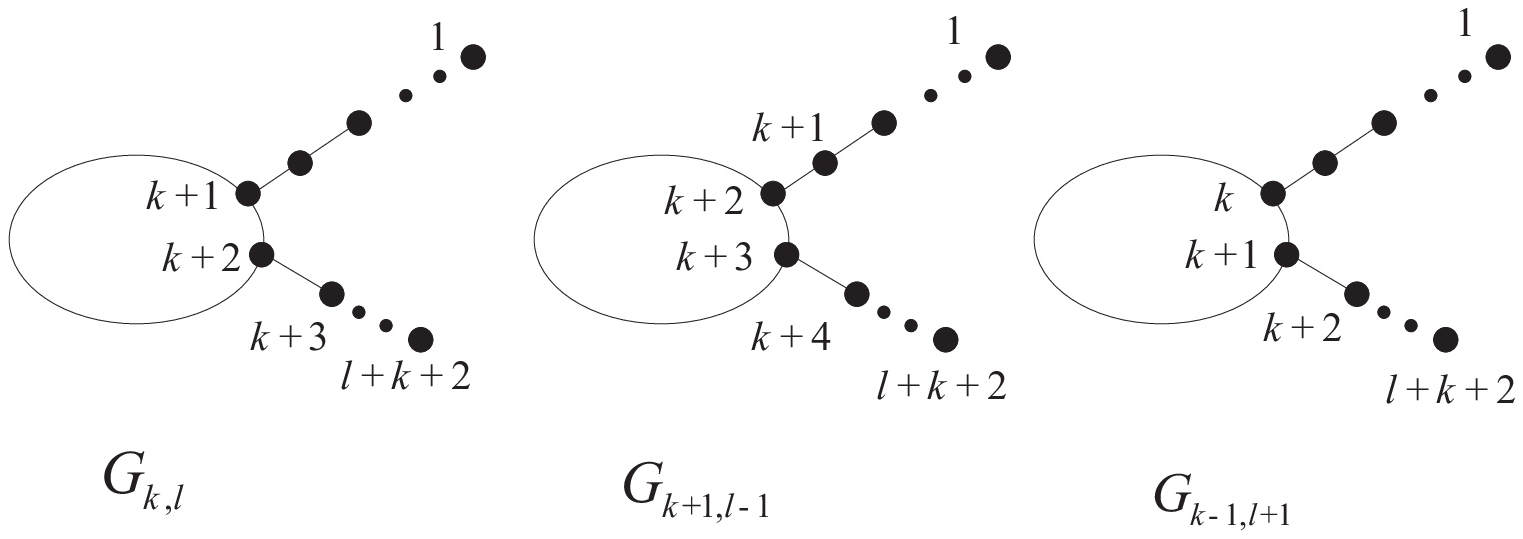}\\
    {$~~$\par\large{Fig. 1. Graphs $G_{k, l}$, $G_{k+1, l-1}$, $G_{k-1,l+1}$. }}
\end{center}
\end{figure}

\noindent \textbf{Proof:} Label vertices of $G_{k, l}$, $G_{k+1,l-1}$, $G_{k-1, l+1}$
as in Figure~1. We partition
$V(G_{k,l})$ into $A \cup \{k+2\}\cup B \cup C \cup D \cup E$, where
\begin{eqnarray*}
    A&=&\{1,\ldots, k+1\},\\
    B&=&\{k+3, \ldots, l+k+2\},\\
    C&=&\{j\mid \dist(j,k+1)+1=\dist(j,k+2)\},\\
    D&=&\{j\mid \dist(j,k+1)=\dist(j,k+2)\}, \\
    E&=&\{j\mid \dist(j,k+1)=\dist(j,k+2)+1\}.
\end{eqnarray*}
When we pass from $G$ to $G'$,  the distances within $A\cup\{k+2\}\cup B$ and
within $C\cup D \cup E$ are unchanged; the distances between $A$ and
$C\cup D \cup E$ , $\{k+2\}$ and $E$ are increased by 1; the
distances between $B$ and $C\cup D \cup E$, $\{k+2\}$ and $C$ are
decreased by 1; the distances between $\{k+2\}$ and  $D$  are
unchanged. If the distance matrices are partitioned according to
$A$, $\{k+2\}$, $B$, $C$, $D$, $E$, their difference is

\noindent\[ D(G_{k+1,l-1})-D(G_{k,l})=\\ \left[
\begin{array}{cccccc}
 0 & 0 & 0 & e_A(e_C)^t  &  e_A(e_D)^t &  e_A(e_E)^t\\
 0 & 0 & 0 & -(e_C)^t &  0 & (e_E)^t \\
0 & 0 & 0 & -e_C(e_C)^t  &  -e_C(e_D)^t &  -e_C(e_E)^t\\
(e_C)^te_A&  -e_C & -(e_C)^te_C & 0& 0& 0\\
(e_D)^t e_A & 0 & -(e_D)^te_C & 0& 0&0\\
(e_E)^te_A & e_E & -(e_E)^te_C & 0 & 0 &0
\end{array}
\right],
\]
where
\[
    e_i=(\underbrace{1, \ldots, 1}\limits_{|i|})^t
\]
and  $i=A, C, D, E$. Let $x=(x_1, \ldots, x_n)^t$ be a positive
eigenvector corresponding to $\Lambda_1(G_{k,l})$. Then we have
\begin{eqnarray}
\frac{1}{2}\big{(}\Lambda_1(G_{k+1,l-1})-\Lambda_1(G_{k,l})\big{)}
&\geqslant&\frac{1}{2}x^t\big{(}D(G_{k+1,l-1})-D(G_{k,l})\big{)}x\\\nonumber
&=&\left(\sum\limits_{j\in A}x_j-\sum\limits_{j\in
B}x_j\right)\sum\limits_{j\in C\cup D\cup E}x_j\\\nonumber&
&-x_{k+2}\sum\limits_{j\in C}x_j +x_{k+2}\sum\limits_{j\in
E}x_j\\\nonumber &\geqslant&\left(\sum\limits_{j\in
A}x_j-\sum\limits_{j\in B}x_j-x_{k+2}\right)\sum\limits_{j\in C\cup
D\cup E}x_j\label{21}.
\end{eqnarray}

Similarly, we have
\begin{eqnarray*}
\frac{1}{2}\big{(}\Lambda_1(G_{k-1,l+1})-\Lambda_1(G_{k,l})\big{)} &\geqslant&\frac{1}{2}x^t\big{(}D(G_{k-1,l+1}))-D(G_{k,l})\big{)}x\\
&=&\left(\sum\limits_{j\in B}x_j+x_{k+2}-\sum\limits_{j\in
A}x_j+x_{k+1}\right)\sum\limits_{j\in C\cup D\cup E}x_j\\&
&+x_{k+1}\sum\limits_{j\in C}x_j
-x_{k-2}\sum\limits_{j\in E}x_j\\
&\geqslant&\left(\sum\limits_{j\in B}x_j+x_{k+2}-\sum\limits_{j\in
A}x_j\right)\sum\limits_{j\in C\cup D\cup E}x_j.
\end{eqnarray*}

Since $D(G_{k,l})x=\Lambda_1(G_{k,l})x$, i.e.,
$\Lambda_1(G_{k,l})x_i=\sum \limits_{j=1}^{n}\dist(i, j)x_j$, for
$1\leqslant i \leqslant n$,  we get:

\noindent if $1\leqslant i \leqslant l+1$, then
\begin{eqnarray}
\Lambda_1(G_{k,l})(x_i-x_{l+k+3-i})&=&\sum
\limits_{m=1}^{s}(l+k-2m+3)(x_{l+k+3-m}-x_{m})+(k-l)\sum\limits_{j\in
D}x_j \\\nonumber & &+(k-l-1)\sum \limits_{j\in
C}x_j+(k-l+1)\sum\limits_{j\in E}x_j
\\\nonumber & &-2\sum \limits_{m=1}^{i-1}(x_{l+k+3-m}-x_{m});\label{111}
\end{eqnarray}

\noindent if $l+2\leqslant i \leqslant s$, then
\begin{eqnarray}
\Lambda_1(G_{k,l})(x_i-x_{l+k+3-i})&=&\sum
\limits_{m=1}^{s}(l+k-2m+3)(x_{l+k+3-m}-x_{m})
\\\nonumber& &+(k+l+3-2i)\sum\limits_{j\in C\cup D \cup E}x_j-2\sum
\limits_{m=1}^{i-1}(x_{l+k+3-m}-x_{m});\label{112}
\end{eqnarray}
where $s=\frac{l+k+2}{2}$, if $l+k$ is even; and
$s=\frac{l+k+1}{2}$, if $l+k$ is odd.

From Eqs. (\ref{111}) and (\ref{112}), we get:

\noindent if $2\leqslant i \leqslant l+1$, then
\begin{eqnarray}
&
&\Lambda_1(G_{k,l})(x_{i-1}-x_{l+k-i+4})-\Lambda_1(G_{k,l})(x_i-x_{l+k+3-i})
\\\nonumber& &=2\sum \limits_{r=1}^{i-1}(x_{l+k+3-r}-x_{r}); \label{24}
\end{eqnarray}
\noindent if $i = l+2$, then
\begin{eqnarray}
&
&\Lambda_1(G_{k,l})(x_{i-1}-x_{l+k-i+4})-\Lambda_1(G_{k,l})(x_i-x_{l+k+3-i})
\\\nonumber& &=2\sum \limits_{r=1}^{i-1}(x_{l+k+3-r}-x_{r})+\sum\limits_{j\in
D}x_j+2\sum\limits_{j\in E}x_j; \label{27}
\end{eqnarray}
\noindent if $l+3\leqslant i \leqslant s$, then
\begin{eqnarray}
&
&\Lambda_1(G_{k,l})(x_{i-1}-x_{l+k-i+4})-\Lambda_1(G_{k,l})(x_i-x_{l+k+3-i})
\\\nonumber & &=2\sum \limits_{r=1}^{i-1}(x_{l+k+3-r}-x_{r})+2\sum\limits_{j\in C\cup D \cup
E}x_j.\label{28}
\end{eqnarray}

Equation (\ref{24}) implies that, for $1\leqslant i\leqslant l+1$,
the differences $x_{l+k+3-i}-x_{i}$ are either all non-positive or
non-negative.

\medbreak \textbf{Case 1.} $k=l$.

 If $x_{l+k+3-i}-x_{i}\leqslant 0$, for all $1\leqslant i\leqslant l+1$, then
$$\sum\limits_{j\in A}x_j-\sum\limits_{j\in B}x_j-x_{k+2}\geqslant
0.$$ So we get $\Lambda_1(G_{k,l})\leqslant \Lambda_1(G_{k+1,l-1})$.
Here, if $\Lambda_1(G_{k,l})= \Lambda_1(G_{k+1,l-1})$, then the last
equality  in (\ref{21}) holds, which implies that
$$2\sum\limits_{j\in E}x_j+\sum\limits_{j\in D}x_j=0,$$ i.e.,
$D=E=\emptyset$. Since
 $i-1 \leqslant l=s-1$, for $1 \leqslant m \leqslant i-1$, $l+k-2m+3\geqslant 2 $,  we get
\begin{eqnarray*}
\sum \limits_{m=1}^{s}(l+k-2m+3)(x_{l+k+3-m}-x_{m})-2\sum
\limits_{m=1}^{i-1}(x_{l+k+3-m}-x_{m})\leqslant 0.
\end{eqnarray*}
Combining this with the fact that  $x$ is a positive eigenvector
corresponding to  $\Lambda_1(G_{k,l})$, we get that the right side
of equation (\ref{111}) is strictly less than 0, which contradicts
the fact that the left side of equation (\ref{111})
$\Lambda_1(G_{k,l})(x_i-x_{l+k+3-i}) \geqslant 0$. So we get
$\Lambda_1(G_{k,l})< \Lambda_1(G_{k+1,l-1})$.

If $x_{l+k+3-i}-x_{i}\geqslant 0$, for all $1\leqslant i\leqslant
l+1$, we can get $\Lambda_1(G_{k,l})< \Lambda_1(G_{k-1,l+1})$
similarly. \medbreak \textbf{Case 2.} $k > l$.

If $x_{l+k+3-i}-x_{i}\geqslant 0$, for all $1\leqslant i\leqslant
l+1$, then  from Eqs. (\ref{27}) and (\ref{28}), we can get that
$x_{l+k+3-i}-x_{i}\geqslant 0$, for all $l+2\leqslant i \leqslant
s$. Similar to the above case, we get that the right side of
equation (\ref{24}) is strict larger than 0, which contradicts the
fact that the left side of equation (\ref{111})
$\Lambda_1(G_{k,l})(x_i-x_{l+k+3-i}) \leqslant 0$. So this case is
impossible.

If $x_{l+k+3-i}-x_{i}< 0$, then  $\sum\limits_{j\in
A}x_j-\sum\limits_{j\in B}x_j-x_{k+2}> 0$, which implies  that
$$\Lambda_1(G_{k,l})> \Lambda_1(G_{k+1,l-1}).$$

 This completes the proof.
\hfill$\Box$

From the proof of above theorem, we get the following corollary.
\begin{cor}
Let $v_l$ and $v_m$ be two adjacent vertices of connected graph $G$.
Let $P_l$ and $P_m$ be two pendant paths with roots $v_l$ and $v_m$,
respectively. If $l > m$, then
$$\sum\limits_{j\in V(P_l)} x_j > \sum\limits_{j\in
V(P_m)} x_j.$$
\end{cor}

\begin{thm}\label{2}
 Let  $C_1$  be a component of $G-u$ and
$v_1,\ldots,v_k$  $(1\leqslant k\leqslant
 \deg_G(u)-\deg_{C_1}(u))$ be some vertices of $N_{G}(u)\setminus
N_{C_1}(u)$. Suppose $N_{C_1}(u)\setminus \{v\}=N_{C_1}(v)$, where
$v$ is a vertex of $C_1$ adjacent to $u$. Let $G'$ be the graph
obtained from $G$ by deleting the edges $uv_s$ and adding the edges
$vv_s$ $(1\leqslant s \leqslant k)$. If there exists a  vertex $w
\in V(G)\setminus (V(C_1)\cup \{u\})$ such that $\dist_G(w,v_s)<
\dist_{G'}(w,v_s)$, for all $1\leqslant s \leqslant k$, then
$\Lambda_1(G)<\Lambda_1(G')$.
\end{thm}

\noindent \textbf{Proof:}  We partition $V(G)$ into $A\cup B \cup
\{u\}\cup \{v\}\cup C$, where
\begin{eqnarray*}
A&=&\{v_1,\ldots,v_k\},
\\ B&=&V(G)\setminus (A\cup V(C_1)\cup \{u\}),
\\
C&=&V(C_1)\setminus \{v\}.
\end{eqnarray*}
From $G$ to $G'$, the distances between $A\cup B$ and $C$ are
unchanged; the distances between $A$ and $\{u\}\cup \{w\}$ are
increased by 1; the distances between $A$ and $\{v\}$ are decreased
by 1; the distances between $A$ and $B\setminus \{w\}$ are not
decreased. Let $x=(x_1, \ldots, x_n)^t$ be a positive eigenvector
corresponding to $\Lambda_1(G)$. Similar to the proof of above
theorem, we have
\begin{eqnarray*}
\frac{1}{2}\big{(}\Lambda_1(G')-\Lambda_1(G)\big{)} &\geqslant&\frac{1}{2}x^t\big{(}D(G')-D(G)\big{)}x\\
&\geqslant& (x_w +x_u-x_{v})\sum\limits_{j\in A}x_j.
\end{eqnarray*}
Since $D(G)x=\Lambda_1(G)x$, i.e., $\Lambda_1(G)x_i=\sum
\limits_{j=1}^{n}\dist(i,j)x_j$,  we can easily get
\begin{eqnarray*}
(\Lambda_1(G)+1)(x_w +x_u-x_v)&=&\sum
\limits_{j=1}^{n}\left(\dist(w,j)+\dist(u,j)-\dist(v,j)\right)x_j+(x_w
+x_u-x_v)\\ &=&\sum \limits_{j\in A\cup (B\setminus
\{w\})}(\dist(w,j)+\dist(u,j)-\dist(v,j))x_j\\
& &+ \sum \limits_{j\in
C}(\dist(w,j)+\dist(u,j)-\dist(v,j))x_j\\
& &+(\dist(w,w)+\dist(u,w)-\dist(v,w))x_w\\
& &+(\dist(w,u) +\dist(u,u)-\dist(v,u))x_u\\&
&+(\dist(w,v)+\dist(u,v)-\dist(v,v))x_v+(x_w +x_u-x_v).
\end{eqnarray*}
As we know, $\dist(u,j)-\dist(v,j)=-1$ and $\dist(w,j)\geqslant 1$,
for $j\in A\cup (B\setminus \{w\})$, so $$\sum \limits_{j\in A\cup
(B\setminus \{w\})}(\dist(w,j)+\dist(u,j)-\dist(v,j))x_j \geqslant
0.$$ Similarly, we can get $$\sum \limits_{j\in
C}(\dist(w,j)+\dist(u,j)-\dist(v,j))x_j>0.$$ Since
\begin{eqnarray*}
\dist(w,w)+\dist(u,w)-\dist(v,w)&=&-1,\\
\dist(w,u) +\dist(u,u)-\dist(v,u)&\geqslant& 0,\\
\dist(w,v)+\dist(u,v)-\dist(v,v)&\geqslant& 3,
\end{eqnarray*}
combining these with all the above equations and inequations,  we
get
$$\big{(}\Lambda_1(G)+1)(x_w +x_u-x_v\big{)}> 0,$$

\noindent which implies $\Lambda_1(G')> \Lambda_1(G)$.

This completes the proof. \hfill$\Box$

\section{Applications}

The graph $G_{n,k}$  is a graph obtained by adding paths $P_{l_1+1},
\ldots, P_{l_{n-k}+1}$ of almost equal lengths (by the length of a
path, we mean the number of its vertices) to the vertices of the
complete graph $K_{n-k}$; that is, the lengths $l_1,\ldots, l_{n-k}$
of $P_{l_1+1}, \ldots, P_{l_{n-k}+1}$ which satisfy
$|l_i-l_j|\leqslant 1$; $1\leqslant i, j\leqslant n-k$.

$K_{n}^k$ is a graph obtained by joining $k$ independant vertices to
one vertex of $K_{n-k}$.

\begin{thm}\label{5}
Of all the connected graphs with $n$ vertices and $k$ cut vertices,
the minimal distance spectral radius is obtained uniquely at
$G_{n,k}$.
\end{thm}
\noindent \textbf{Proof:}  We are supposed to prove that if $G$ is a
connected graph with $n$ vertices and $k$ cut vertices, then
$\Lambda_1(G)\geqslant \Lambda_1(G_{n,k})$ with equality only when
$G\cong G_{n,k}$. Let $V_1$ be the set of the cut vertices of $G$.
Note that if we add some edges to $G$ such that each block of
$G-V_1$ is a clique, denoting the new graph by $G'$, then
$\dist_G(i,j)\geqslant \dist_{G'}(i,j)$. Consequently,
$D(G)\geqslant D(G')$, which implies $\Lambda_1(D(G))\geqslant
\Lambda_1(D(G'))$ with equality only when $D(G)=D(G')$. So, in the
following, we always assume that each cut vertex of $G$ connects
exactly two blocks and that all these blocks are cliques.
 Order the
cardinalities of these blocks $a_1\geqslant a_2\geqslant \cdots
\geqslant a_{k+1}\geqslant 2$ and denote the blocks by  $K_{a_1},
\ldots, K_{a_{k+1}}$.  If $k=0$, then $G\cong K_{n}$ and the theorem
holds. If $k=n-2$, then $G$ is the path $G_{n,n-2}$. If $k=n-3$,
then $a_1=3$, $a_2= \cdots = a_{k+1}=2$. The result follows from a
repeated use of Theorem \ref{1}. Thus we may assume that $1\leqslant
k\leqslant n-4$. Moreover, we observe that
$a_1=n+k-(a_2+\cdots+a_{k+1})\leqslant n-k$.

Choose $G$ such that the distance spectral radius is as small as
possible.

\medbreak \textbf{Claim.} $a_1=n-k$.

Otherwise, $a_1\leqslant n-k-1$, which implies $a_2\geqslant \cdots
\geqslant a_i\geqslant 3$ for some $i$, $2\leqslant i\leqslant k+1$.

Suppose $K_{a_{i_1}}, \ldots, K_{a_{i_t}}$ are the blocks, each of
which contains  at least two roots of pendant paths of $G$. Let
$\mathcal{P}$ be the set of pendant paths whose roots are contained
in $K_{a_{i_1}}, \ldots, K_{a_{i_t}}$ and $P_{m}$ be one of the
shortest pendant paths among $\mathcal{P}$.  Suppose the root of
$P_{m}$ is contained in $K_{a_{i_s}}$ for some $1\leqslant s
\leqslant t$ and $P_{l}$ is another pendant path whose root is also
contained in $K_{a_{i_s}}$. Then we have $0\leqslant l-m \leqslant
1$. Otherwise, by Theorem \ref{1}, we can find a graph $G'$ such
that $\Lambda_1(G)> \Lambda_1(G')$, which is a contradiction.  We
label the vertices of $G$ such that
\begin{eqnarray*}
P_m&=&v_1\cdots v_{m},\\ P_l&=&v_{m+1}\cdots v_{m+l},\\
V(K_{a_{i_s}})\setminus\{v_{m}, v_{m+1}\}&=&\{u_1, \ldots,
u_r\},
\end{eqnarray*}
where $v_{m}$ and $v_{m+1}$ are the roots of $P_m$ and $P_l$,
respectively. Suppose
 $u_{1}$ is a cut vertex
of $G$ such that $G[V(C_1)\cup \{u_1\}]$ contains at least one block
$K_{a_h}$, for some $1 \leqslant h \leqslant i$, where $C_1$ is the
component of $G-u_1$ which does not contain $v_m$.

Let $N_{C_1}(u_1)=\{w_1. \ldots, w_q\}$ and
\begin{eqnarray*}
G'=G&-&v_{m}u_2-\cdots
-v_{m}u_{r}-v_{m}v_{m+1}\\
&
+&w_1u_2+\cdots+w_1u_{r}+w_1v_{m+1}+\\
& &\cdots\cdots\cdots\\&
+&w_{q-1}u_2+\cdots+w_{q-1}u_{r}+w_{q-1}v_{m+1}\\
&+&w_qu_2+\cdots+w_qu_{r}+w_qv_{m+1}.
\end{eqnarray*}
Then $G'$ is a connected graph with $n$ vertices and $k$ cut
vertices.

In the following, we consider the graph $G'$. Since $V(G')=V(G)$, we
partition $V(G')$ into $A \cup \{u_{1}\} \cup B \cup C $, where
\begin{eqnarray*}
A&=&V(P_m),\\ B&=&V(C_1),\\ C&=&V(G)\setminus (A
\cup \{u_{1}\} \cup B).
\end{eqnarray*}
 From $G'$ to $G$, the distances between $A$ and
$B\cup \{u_{1}\}$, $\{u_{1}\}$ and $B$ are unchanged; the distances
between $B$ and $C$ are increased by 1; the distances between $A$
and $C$ are decreased by 1. Let $x=(x_1, \ldots, x_n)^t$ be a
positive eigenvector corresponding to $\Lambda_1(G')$. Then we have
\begin{eqnarray}
\frac{1}{2}\big{(}\Lambda_1(G)-\Lambda_1(G')\big{)}
&\geqslant&\frac{1}{2}x^t\big{(}D(G)-D(G')\big{)}x\\\nonumber
&\geqslant& \left(\sum\limits_{j\in B}x_j-\sum\limits_{j\in
A}x_j\right)\sum\limits_{j\in C}x_j.\label{117}
\end{eqnarray}

\textbf{Case 1.} $\sum\limits_{j\in B}x_j-\sum\limits_{j\in A}x_j >
0$.

From inequality (\ref{117}), we get that $\Lambda_1(G)>
\Lambda_1(G')$, which is a contradiction.

\medbreak \textbf{Case 2.} $\sum\limits_{j\in
B}x_j-\sum\limits_{j\in A}x_j \leqslant 0$

 Since $D(G')x=\Lambda_1(G')x$, i.e., $\Lambda_1(G')x_i=\sum
\limits_{j=1}^{n}\dist(i,j)x_j$,  we can easily get that
\begin{eqnarray}
\Lambda_1(G')\left(\sum\limits_{j\in B}x_j-\sum\limits_{j\in
A}x_j\right) &=& \sum \limits_{i\in A\cup \{u_1\}}\left(\sum
\limits_{u\in B}\dist(u, i)- \sum \limits_{v\in
A}\dist(v,i)\right)x_i\\\nonumber & &+\sum \limits_{i\in
B}\left(\sum \limits_{u\in B}\dist(u,i)- \sum \limits_{v\in
A}\dist(v,i)\right)x_i\\\nonumber & &+\sum \limits_{i\in
C}\left(\sum \limits_{u\in B}\dist(u, i)- \sum \limits_{v\in
A}\dist(v,i)\right)x_i.\label{113}
\end{eqnarray}
Since $G[V(C_1)\cup \{u_1\}]=G'[V(C_1)\cup \{u_1\}]$ contains at
least one block $K_{a_h}$, for some $1 \leqslant h \leqslant i$,
$G'[V(C_1)\cup \{u_1\}]$ must contain at least two pendant paths
$P'$ and $P''$  whose roots are contained in the same block. Denote
the roots of $P'$ and $P''$ by $\omega_1$ and $\omega_2$,
respectively. Suppose $dist_{G'}(v_{m+1}, \omega_1)=k$, then $k
\geqslant 2$ and $dist_{G'}(v_{m+1}, \omega_2)=k$. As we know, $P'$
and $P''$ are two pendant paths of length at least $m$, so

\noindent if $i\in A\cup \{u_1\}$, then
\begin{eqnarray}\sum \limits_{u\in B}\dist(u, i)-
\sum \limits_{v\in A}\dist(v,i) \geqslant 2 \times
\frac{(m+k-1)(m+k)}{2}-\frac{m(m-1)}{2};\label{114}\end{eqnarray} if
$i\in B$, then
\begin{eqnarray}
\sum \limits_{u\in B}\dist(u, i)- \sum \limits_{v\in
A}\dist(v,i)&\geqslant& \frac{(k-2)(k-1)}{2}-(m-1)k ;\label{115}
\end{eqnarray}
 if $i\in C$, then
\begin{eqnarray}\sum \limits_{u\in B}\dist(u,i)- \sum \limits_{v\in
A}\dist(v,i)> 0.\label{116}\end{eqnarray} is obvious.

 Combining Eqs.
(\ref{113}), (\ref{114}), (\ref{115}) with  (\ref{116}), we get
\begin{eqnarray*}
 & &\Lambda_1(G')\left(\sum\limits_{j\in B}x_j-\sum\limits_{j\in A}x_j\right)\\
& & > \left[2 \times
\frac{(m+k-1)(m+k)}{2}-\frac{m(m-1)}{2}\right]\sum \limits_{i\in
A\cup \{u_1\}}x_i\\&
&~~~+\left[\frac{(k-2)(k-1)}{2}-(m-1)k\right]\sum \limits_{i\in
B}x_i\\& &> \left[
(m+k-1)(m+k)-\frac{m(m-1)}{2}+\frac{(k-2)(k-1)}{2}-(m-1)k\right]\sum
\limits_{i\in B}x_i\\& &> 0,
\end{eqnarray*}
 which is a contradiction. So this case does not exist.

 Up to now, we have proved the claim that  $a_1=n-k$, which implies
 that $$a_2= \cdots = a_{k+1}=2.$$ From a repeated use of Theorem
 \ref{1}, we get that $G\cong
G_{n,k}$.

This completes the proof. \hfill$\Box$

\begin{thm}
Of all the connected graphs with $n$ $(n\geqslant 4)$ vertices and
$k$ cut edges, the minimal distance spectral radius is obtained
uniquely at $K_{n}^k$.
\end{thm}
\noindent \textbf{Proof:} We are supposed to prove that if $G$ is a
connected graph with $n$ vertices and $k$ cut edges, then
$\Lambda_1(G)\geqslant \Lambda_1(K_{n}^k)$ with equality only when
$G\cong K_{n}^k$. Let $E_1 = \{e_1, e_2, \ldots, e_k\}$ be the set
of the cut edges of $G$. For a similar reason to Theorem \ref{5}, we
assume that each component of $G-E_1$ is a clique.

If $k=0$, then $G\cong K_{n}$ and the theorem holds. So we assume
that $k\geqslant 1$. Denote the components of $G-E_1$ by $K_{a_0},
\ldots, K_{a_k}$, where $a_0+\cdots +a_k=n$.

Let $V_{a_i}$ = $\{v\in K_{a_i}:$ $ v$ is an end vertex of the  cut
edges of  $G\}$, and choose $G$ such that the distance spectral
radius is as small as possible.

\medbreak\textbf{Claim 1.} $|V_{a_i}|=1$, $0\leqslant i \leqslant
k$.

Otherwise, $|V_{a_i}|>1$ for some  $0\leqslant i \leqslant k$.
Suppose $u$, $v \in V_{a_i}$ and $vv_j, uv_h\in E_1$. Let
$$G'=G-vv_j+uv_j.$$
Then, $G'$ is still a connected graph with $n$ vertices and $k$ cut
edges, and $\dist_{G'}(v_h, v_{j}) < \dist_G(v_h, v_{j})$.  In $G'$,
$N_{K_{a_i}}(u)\setminus\{v\}=N_{K_{a_i}}(v)\setminus\{u\}$, so by
Theorem \ref{2}, we get $\Lambda_1(G')< \Lambda_1(G)$, which is a
contradiction.

So, in the following, we can assume that $V_{a_i}=\{v_i\}$,
$0\leqslant i \leqslant k$.

\medbreak\textbf{Claim 2.} If $v_s$ is adjacent to $v_t$, where
$0\leqslant t, s \leqslant k$, then $\dg(v_s)=1$ or $\dg(v_t)=1$.

Otherwise, $\dg(v_s)\geqslant 2$ and $\dg(v_t)\geqslant 2$.  Denote by
$C_1$, $C_2$ the  components of $G-v_sv_t$ which contain $v_s$ and
$v_t$, respectively.    Suppose $N(v_s)\setminus \{v_t\}=\{v_{k+1},
\ldots, v_{k+r}\}$ and $u\in N(v_t)\setminus \{v_s\}$. Let
$$G'=G-v_sv_{k+1}-\ldots-v_sv_{k+r}+v_tv_{k+1}+\ldots+v_tv_{k+r}.$$
Then, $G'$ is still a connected graph with $n$ vertices and $k$ cut
edges, and $\dist_{G'}(u, v_{k+i}) < \dist_G(u, v_{k+i})$ $(1
\leqslant i \leqslant r)$. In $G'$, $v_sv_t$ is a pendant edge, so
by Theorem \ref{2}, we get $\Lambda_1(G')< \Lambda_1(G)$, which is a
contradiction.

Since $G$ is connected, combining Claim 1 with Claim 2, we get that
$G\cong K_{n}^k$.

This completes the proof. \hfill$\Box$


\end{document}